\documentclass[12pt,leqeq]{article}

\usepackage{amssymb}

\usepackage{amsmath}

\usepackage[all]{xy}

\newtheorem{theorem}{Theorem}

\newtheorem{corollary}[theorem]{Corollary}
\newtheorem{definition}[theorem]{Definition}
\newtheorem{example}[theorem]{Example}

\numberwithin{equation}{section}

 \newcommand{\kr}{\mathbf{k}}
 \newcommand{\Z}{\mathbb{Z}}
 
 \newcommand{\Agr}{A^\bullet}

 \newcommand{\E}{\mathcal{E}}
 
 \newcommand{\N}{\mathbb{N}}
 
 \newcommand{\R}{\mathbb{R}}
 \newcommand{\C}{\mathbb{C}}

\begin{document}

\setlength{\baselineskip}{16pt}

\title{N-flat connections}

\author{Mauricio Angel\hspace{.3cm} and Rafael D\'\i az}

\maketitle

%%% ----------------------------------------------------------------------
\begin{sloppypar}

\begin{abstract}
We construct geometric examples of N-differential graded algebras
such as the algebra of differential forms of depth $N$ on an
affine manifold, and $N$-flat covariant derivatives.
\end{abstract}

%%% ----------------------------------------------------------------------

\section*{Introduction}

The applications of the theory of complexes and homological
algebra have touched many branches of mathematics such as
topology, geometry and mathematical physics. The possibility of
developing an homological algebra for the equation $d^N=0$ for
$N\geq 3$ has been around at least since 1940 in works by Mayer
[M].

\smallskip

However, N-homological algebra for $N\geq 3$ only aroused the
interest that it deserves in 1991 when Kapranov's paper [K]
appeared. Fundamental works by Dubois-Violette appeared in [DV1],
[DV2] soon after. The key notion is that of a q-differential
graded algebra (q must be a primitive N-root of unity), which
consists of a $\Z$-graded vector space $V$ together with an
operation $m:V\otimes V\to V$ and an degree one map $d:V\to V$
such that: 1) $m$ is associative, 2) $d$ satisfy the q-Leibniz
rule $d(ab)=d(a)b+q^{deg(a)}ad(b)$ and 3) $d^N=0$.

\smallskip

Notice that in the definition of a q-differential graded algebra
not only the equation $d^2=0$ is deformed to $d^N=0$, but also the
Leibniz rule is substituted by the q-Leibniz rule. The theory of
q-differential graded algebras has been further developed in many
papers such as [AB], [DV3], [DVK], [KW], [KN], [S].

\smallskip

The foundations of a theory of N-differential graded algebras
(N-dga) has been written down in [AD]. A N-dga consists of a
$\Z$-graded vector space $V$ together with an associative
operation $m:V\otimes V\to V$ and a degree one map $d:V\to V$ such
that: 1) $d^N=0$ and 2) $d(ab)=d(a)b+(-1)^{deg(a)}ad(b)$. Notice
that in the definition of a N-dga the equation $d^2=0$ is replaced
by $d^N=0$, whereas the Leibniz rule is not modified. This fact
explain why our definition is better suited for the use of
differential geometric techniques.

\smallskip

At first sight it looks difficult to came up with examples of
N-dga. The purpose of this paper is to introduce geometric
examples of N-differential graded algebras, which show that that
this sort of algebras appear naturally in a wide variety of
contexts. Our examples range from the algebra of differential
forms of depth N (see Section 4) to the theory of N-flat covariant
derivatives (see Section 1).

\section{N-covariant differentials}

The purpose of this paper is to introduce geometric examples of
N-differential graded algebras (N-dga). Let us first formally
introduce the notion of N-dga [AD]. Let $\kr$ be a commutative
ring with unit.

\begin{definition}\label{Ndga}
A {\bf N-differential graded algebra} or N-dga over $\kr$, is a
triple $(\Agr,m,d)$ where $m:A^k\otimes A^l\to A^{k+l}$ and
$d:A^k\to A^{k+1}$ are $\kr$-modules homomorphisms satisfying
\begin{enumerate}
\item[1)] $(\Agr,m)$ is a graded associative algebra.

\item[2)] $d$ satisfies the graded Leibniz rule
$d(ab)=d(a)b+(-1)^{\bar{a}}ad(b)$.

\item[3)] $d^{N}=0$, i.e., $(\Agr,d)$ is a $N$-complex.
\end{enumerate}

\end{definition}

A 1-dga is a graded associative algebra. A 2-dga is a differential
graded algebra.

\smallskip

Definition above may be justified categorically as follows:
consider the category ${\mathbf{NComp_\kr}}$ of nilpotent
differential graded $\kr$-modules. Objects in
${\mathbf{NComp_\kr}}$ are pairs $(V,d)$ where $V$ is a
$\Z$-graded $\kr$-module $V=\displaystyle{\oplus_{i\in\Z}V^i}$
together with a degree one map $d:V\to V$ such that $d^N=0$ for
some $N\geq 2$. ${\mathbf{NComp_\kr}}$ is a symmetric monoidal
category since if $(V_1,d_1)$ is such that $d_1^{N_1}=0$ and
$(V_2,d_2)$ is such that $d_2^{N_2}=0$ then $(V_1\otimes
V_2,d_1\otimes Id+Id\otimes d_2)$ satisfy $(d_1\otimes
Id+Id\otimes d_2)^{N_1+N_2-1}=0$. N-differential graded algebras
(for $N\geq 2$) are the monoids in ${\mathbf{NComp_\kr}}$.

\begin{definition}
A N-dga $(\Agr,m,d)$ is called {\bf proper} if $d^{N-1}\neq 0$.
\end{definition}

Let $M$ be a smooth finite dimensional manifold and $\E=M\times E$
be a trivial bundle over $M$. The reader may assume that $\E$ is a
global bundle and that $M$ is actually a neighborhood on which
$\E$ is trivial. Since our results are covariant they will hold
globally as well.

\smallskip

The space $\Omega^{\bullet}(M,End(\E))$ of $End(\E)$-valued forms
on $M$ is endowed with a differential graded algebra structure,
with the product given by
\[\begin{array}{c} \wedge:\Omega^{\bullet}(M,End(\E))\otimes\Omega^{\bullet}(M,End(\E))\to\Omega^{\bullet}(M,End(\E))\\
(\alpha\otimes\psi)(\beta\otimes\phi)\longmapsto
(\alpha\wedge\beta)\otimes \psi\phi.\\ \end{array}\] The space of
$\E$-valued forms $\Omega^{\bullet}(M,\E)$ is endowed with a
differential graded module structure over
$\Omega^{\bullet}(M,End(\E))$ and also over the differential
graded algebra $\Omega^{\bullet}(M)$ of differential forms on $M$,
given by
\[\begin{array}{c} \Omega^{\bullet}(M)\otimes\Omega^{\bullet}(M,\E)\to\Omega^{\bullet}(M,\E)\\
\alpha(\beta\otimes\phi)\longmapsto
(\alpha\wedge\beta)\otimes\phi.\end{array}\]

\[\begin{array}{c}
\Omega^{\bullet}(M,End(\E))\otimes\Omega^{\bullet}(M,\E)\to\Omega^{\bullet}(M,\E)\\
(\alpha\otimes\Phi)(\beta\otimes\phi)\longmapsto
(\alpha\wedge\beta)\otimes\Phi(\phi)\\
\end{array}\]

Recall that a covariant derivative $\nabla$ on $\E$ is a linear
map $\nabla:\Omega^\bullet(M,E)\to\Omega^\bullet(M,E)$  of degree
1 such that
\[\nabla(a\alpha)=(da)\alpha+(-1)^{deg(a)}a\nabla\alpha, \text{for all $a\in\Omega^\bullet(M)$, $\alpha\in\Omega^\bullet(M,E)$}.\]
Since $\E$ is trivial, $\nabla$ may be written as
$\nabla=d+\omega$ for some $\omega\in\Omega^1(M,End(\E))$, where
$d$ is the de Rham differential and $\omega$ is the connection
one-form of $\nabla$.

\smallskip

For any $\alpha\in\Omega^k(M,E)$,
$\nabla\alpha\in\Omega^{k+1}(M,E)$ is the $E$-valued form given by
\begin{equation}
\nabla\alpha=d\alpha+\omega\wedge\alpha.\label{Cov}
\end{equation}
Taking the covariant derivative of (\ref{Cov}) we get
\begin{equation}
\nabla^2(\alpha)=d\omega\wedge\alpha+\omega\wedge\omega\wedge\alpha=(d\omega+\omega\wedge\omega)\wedge\alpha=F_\omega\wedge\alpha.
\end{equation}
The 2-form $F_\omega\in\Omega^2(M,End(\E))$ is called the
curvature of $\nabla$. A connection $\omega$ is said to be {\bf
flat} if $\nabla^2=0$ or equivalently if $F_\omega=0$.

\begin{definition}
We say that a connection $\omega$ is {\bf N-flat} if
$\nabla^{N}=0$, where $\nabla=d+\omega$ and $N\geq 2$.
\end{definition}

A 2-flat connection is just a flat connection. We have the
following
\begin{theorem}\label{MC2N}
Let $M$ be a manifold and $\mathcal{E}=M\times E$ be a trivial
bundle over $M$. A connection $\omega$ is
\begin{enumerate}
\item $2N$-flat if and only if $F_\omega^N=0$. \item $(2N+1)$-flat
if and only if $F_\omega^N\nabla=0$.
\end{enumerate}
\end{theorem}
{\bf proof} Suppose that $K=2N+n$, $n\in\{0,1\}$, then
\[\nabla^{K}=(d+\omega)^{2N+n}=(d(\omega)+\omega\wedge\omega)^N\nabla^n=F_\omega^N\nabla^n.\blacklozenge\]

\begin{example}
Consider $\R^4$ with coordinates $(x_1,x_2,x_3,x_4)$. For
$\omega_1,\omega_2:\R^4\to \R$ we consider the connection
$\omega=\omega_1dx_1+\omega_2dx_2$.

\smallskip

A simple calculation shows that
\[F_\omega=\Bigl(\frac{\partial\omega_1}{\partial x_2}-
\frac{\partial\omega_2}{\partial x_1}\Bigr)dx_1\wedge dx_2\neq 0,
\ \text{if}\ \frac{\partial\omega_1}{\partial
x_2}-\frac{\partial\omega_2}{\partial x_1}\neq 0.\] For example
one can take $\omega_1=x_2$ and $\omega_2=-x_1$, then
$F_\omega=2dx_1\wedge dx_2\neq 0$. But $F_\omega^2=0$ which
implies that $\omega$ is a 4-flat connection and
$(\Omega^\bullet(\R^4,End(\R^4)),\nabla)$ is a proper 4-dga.
\end{example}

The previous example is an instance of the following result.

\begin{theorem}
Let $M$ be a manifold and assume that $TM=A\oplus B$, where
$dim(A)=a$, $dim(B)=b$, and $\omega$ is a connection on $TM$ such
that if $F_\omega(\alpha,\beta)\neq 0$ then $\alpha\in A$ and
$\beta\in A$. In this case  $\nabla=d+\omega$ is $2N$-flat for
$N>a$.
\end{theorem}

Let now $M$ be a $n$-dimensional smooth manifold and $\E=M\times
E$ a trivial bundle on it. Using local coordinates
$x_1,\cdots,x_n$ a connection $\omega\in\Omega^1(M,End(\E))$ may
be written as $\omega=\omega_idx^i$. The 2-form
$d\omega+\omega\wedge\omega$ can be written as $F_{ij}dx^i\wedge
dx^j$ where $F_{ij}=\partial_i \omega_j-\partial_j \omega_i+
[\omega_i,\omega_j]$. Furthermore
\[(F_{ij}dx^i\wedge dx^j)^k=\sum_{A\subseteq [n]}\biggl(\sum_{\alpha\in
P(A)}\prod_{i=1}^k
sign(\alpha)F_{a_i,b_i}\biggr)dx_{s_1}\wedge...\wedge
dx_{s_{2k}}\] where $[n]$ denotes the set $\{1,2,...,n\}$, the
cardinality of $A$ is $2k$ and $P(A)$ is the set of ordered
pairings of $A=\{s_1<\cdots<s_{2k}\}$. An ordered pairing
$\alpha\in P(A)$ is a sequence $\{(a_i,b_i)\}_{i=1}^k$ such that
$A=\bigsqcup_{i=1}^k \{a_i,b_i\}$ and $a_i<b_i$. Theorem
\ref{MC2N} part (1) implies

\begin{theorem}
$(\Omega^{\bullet}(M,\E),d+\omega)$ is a $2k$-dga if and only if
\[\sum_{\alpha\in P(A)}sign(\alpha)\prod_{i=1}^k F_{a_i,b_i}=0,\]
for all $A\subseteq [n]$ with cardinality $2k$.
\end{theorem}

\begin{corollary}
If $dim(M)=2n$, then $(\Omega^{\bullet}(M,\E),d+\omega)$ is a
$2n$-dga if and only if
\[\sum_{\alpha\in P([2n])}sign(\alpha)\prod_{i=1}^n F_{a_i,b_i}=0.\]
\end{corollary}

Let $(M,g)$ be a Riemannian manifold. The tangent bundle $TM$ has
a canonical covariant derivative $\nabla$, called the Levi-Civita
connection. Suppose that the Riemannian metric is given in local
coordinates $x_1,\cdots,x_n$ by the positive definite symmetric
matrix $g_{ij}$:
\[ds^2=\sum g_{ij}dx^idx^j.\]
Denote by $\partial_1,\cdots,\partial_n$ the corresponding
coordinate vector fields, then the covariant derivatives can be
expressed
\[\nabla_{\partial_i}\partial_j=\sum_k\Gamma^i_{jk}\partial_k\hspace{1cm}\Gamma^i_{jk}=\sum_l\frac{1}{2}(\partial_kg_{lj}+\partial_jg_{lk}-\partial_{l}g_{jk})g^{il},\]
where $(g^{lk})$ is the inverse matrix to $(g_{kl})$. In this case
the curvature 2-form is given by
\[R=R_{kl}dx_k\wedge dx_l=R_{jkl}^i\partial_j\otimes dx_i\otimes dx_k\wedge dx_l,\]
where
\[R_{jkl}^i=\partial_k\Gamma^i_{jl}-\partial_l\Gamma^i_{jk}+\Gamma^h_{jl}\Gamma^i_{hk}-\Gamma^h_{jk}\Gamma^i_{hl}.\]
\begin{theorem}
$(\Omega^\bullet(M,TM),\nabla)$ is a $2k$-dga if and only if
\[\sum_{\alpha\in P(A)}sign(\alpha)\prod_{i=1}^k R_{a_ib_ik}^l=0,\]
for all $A\subseteq [n]$ with cardinality $2k$.
\end{theorem}

%\begin{definition}
%We say that a Riemannian manifold if {\bf locally N-flat} if the
%Levi-Civita connection $\nabla$ satisfy $\nabla^N=0$. In the
%2-flat case this notion agrees with the notion of a manifold
%locally flat.
%\end{definition}

\begin{example}
Consider the space $S^2\times T^2$ with coordinates
$(x_1,x_2,x_3,x_4)$, $0< x_1< 2\pi$, $-\pi/2< x_2< \pi/2$,
$0<x_3<1$ and $0<x_4<1$. The metric on $S^2\times T^2$ is given in
those coordinates by
\[ds^2=sin^2(x_2)(dx_1)^2+(dx_2)^2+(dx^3)^2+(dx^4)^2.\]
The curvature is
\[R=\left(\begin{array}{cc} 0& 1 \\ -sin^2(x_2)&0 \\ \end{array}\right)dx_1\wedge dx_2 \neq 0\]
Since $R^2=0$, then $(\Omega(S^2\times\R^2,T(S^2\times
T^2)),\nabla)$ is a 4-dga.
\end{example}

We generalize the previous example as follows:

\begin{theorem}
Let $(M,g)$ be a Riemannian manifold and assume that $TM=A\oplus
B$ is an orthogonal decomposition, where $dim(A)=a$, $dim(B)=b$,
and $g$ is flat on $B$, then $(\Omega(M,TM),\nabla)$ is a $2N$-dga
for $N\geq a$.
\end{theorem}

Recall the definition of the generalized cohomology groups of a
N-complex $(A,d)$ [Kap]
\[_pH^{i}(A)=\frac{Ker\{d^{p}:A^{i}\to A^{i+p}\}}{Im\{d^{N-p}:A^{i-N+p}\to A^{i}\}}.\]
The total object $\mbox{{\bf H$^\bullet$}}(A)$ of the cohomology
associated to a $N$-complex is $\mbox{{\bf
H$^\bullet$}}(A)=\bigoplus_{m=0}^{\infty}\mbox{\bf H$^m$}(A)$,
where
\[\mbox{\bf H$^m$}(A)=\bigoplus_{2i-p=m}{ _pH^{i}(A)}.\]

\begin{example}
Consider the torus $T^4$ with global coordinates
$(\theta_1,\theta_2,\theta_3,\theta_4)$ satisfying
$(\theta_1,\theta_2,\theta_3,\theta_4)=(\theta_1+1,\theta_2+1,\theta_3+1,\theta_4+1)$.
For
\[E_{11}=\Bigl(\begin{matrix} 1&0\\0&0\end{matrix}\Bigr)\hspace{1cm}E_{12}=\Bigl(\begin{matrix} 0&1\\0&0\end{matrix}\Bigr)\]
We define $\nabla=d+E_{11}d\theta_1+E_{12}d\theta_2$, then since
$[E_{11},E_{12}]=E_{12}$ we have $\nabla^2=E_{12}d\theta_1\wedge
d\theta_2\neq 0$, nevertheless it is obvious that $\nabla^4=0$
thus $(\Omega(T^4,\C^2),\nabla)$ is a (proper) 4-dga.

\smallskip

We have $(d\theta_1\wedge d\theta_2)(d\theta_1)=0$ which implies
that $d\theta_1\in Ker(\nabla^2)$, and $d\theta_1$ is not in
$Im(\nabla^{2})$, because otherwise $\nabla^2\beta$ is a 2-form
for all functions $\beta$. Then the cohomology group
$_2H^1(T^4)\neq 0$.

\smallskip

In [AD] we showed that if $(A,d)$ is a N-dga then ({\bf H}$(A),d$)
is a $(N-1)$-dga, thus we see that $\mbox{{\bf H}
}(\Omega(T^4,\C^2),\nabla)(\neq 0)$ is a (proper) 3-dga.
\end{example}

Let $M$ be a manifold and $\E_i$ trivial bundles on $M$, $i=1,2$.
For $\nabla^{\E_i}$ covariant derivatives on $\E_i$, we have a
natural covariant derivative $\nabla^{\E_1\otimes\E_2}$ on the
tensor product $\E_1\otimes\E_2$ given by
\[\nabla^{\E_1\otimes\E_2}(s_1\otimes s_2)=\nabla^{\E_1}s_1\otimes s_2+s_1\otimes\nabla^{\E_2}s_2,\]
for $s_i\in\Omega(M,\E_i)$.

\smallskip

The next result is an analog to Theorem 8 [AD] and we omit the
proof,

\begin{theorem}
If $(\Omega(M,\E_1),\nabla^{\E_1})$ is a N-dga and
$(\Omega(M,\E_2),\nabla^{\E_2})$ is a M-dga, then
$(\Omega(M,\E_1\otimes\E_2),\nabla^{\E_1\otimes\E_2})$ is a
$(M+N-1)$-dga. Equivalently, if $\nabla^{\E_1}$ is $M$-flat and if
$\nabla^{\E_2}$ is $N$-flat then $\nabla^{\E_1\otimes\E_2}$ is
$(M+N-1)$-flat.
\end{theorem}

\section{N-Chern-Simons actions}

Let $M$ be a manifold, $dim(M)=2K+1$ and $\mathcal{E}=M\times E$
be a trivial bundle over $M$. For a connection $\omega$ on $\E$,
consider its K-Curvature
\(F_{\omega}^{K}=(d\omega+\omega\wedge\omega)^{K}\). On
$\Omega^\bullet(M,E)$ exists a linear functional $\int_M
Tr\!\!:\Omega^\bullet(M,End(E))\to \R$ of degree $2K\!+\!1$, i.e.,
$\int b=0$ if $\bar{b}\neq 2K+1$, given by
\[\omega\longmapsto\int_MTr(\omega), \text{for all $\omega\in\Omega^\bullet(M,End(E))$}.\]
The functional $\int_MTr$ satisfy the conditions of [AD],
\begin{enumerate}
\item $\int_MTr$ is non degenerate, that is,
$\int_MTr(\alpha\wedge\beta)=0$ for all $\alpha$, then $\beta=0$.
\item $\int_M Tr(d(\alpha))=0$ for all $\alpha$, where
$d=d_{End(\E)}$. \item $\int_MTr$ is cyclic, this is
\[\int_MTr(\alpha_1\alpha_2...\alpha_n)=(-1)^{\bar{\alpha_1}(\bar{\alpha_2}+\cdots+\bar{\alpha_n})}\int_MTr(
\alpha_2...\alpha_n\alpha_1).\]
\end{enumerate}
We define the {\bf Chern-Simons} functional
$cs_{2,2K}:\Omega^\bullet(M,End(E))\to\R$ by
\[cs_{2,2K}(\omega)=2K\!\int_M Tr(\!\pi (\#^{-1}(\omega(d\omega+\omega^2)^K)))\]
where
\begin{enumerate}
\item $\R\!\!<\!\!\omega,d\omega\!\!>$ denotes the free
$\R$-algebra generated by symbols $\omega$ and $d\omega$. \item
$\#:\R\!<\!\omega,d\omega\!>\longrightarrow\R\!<\!\omega,d\omega\!>$
is the linear map defined by
\[\#(\omega^{i_1}d(\omega)^{j_1}...\omega^{i_k}d(\omega)^{j_k})=
(i_1+..+i_k+j_1+..+j_k)\omega^{i_1}d(\omega)^{j_1}..
\omega^{i_k}d(\omega)^{j_k}.\] \item
$\pi:\R\!<\!\omega,d(\omega)\!>\longrightarrow
\Omega^\bullet(M,End(E))$ is the canonical projection.
\end{enumerate}

We have the following result
\begin{theorem}
Let $K\geq 1$ be an integer. The Chern-Simons functional
$cs_{2,2K}$ is the Lagrangian for the 2K-Maurer-Cartan equation,
i.e., $\omega\in \Omega^\bullet(M,End(E))$ is a critical point of
$cs_{2,2K}$ if and only if $F_\omega^{K}=0$.
\end{theorem}

For $K=2,\ 3,\ 4$ the Chern-Simons functional $cs_{2,2K}(a)$ is
given by
\begin{eqnarray*}
cs_{2,4}(\omega)&=&\int_MTr(\frac{4}{3}\omega(d(\omega))^2+2\omega^3d(\omega)+\frac{4}{5}\omega^5).\\
cs_{2,6}(\omega)&=&\int_MTr(\frac{3}{2}\omega(d(\omega))^3+\frac{18}{5}\omega^3(d(\omega))^2+
3\omega^5d(\omega)+\frac{6}{7}\omega^7). \\
cs_{2,8}(\omega)&=&\int_MTr(\frac{8}{5}\omega(d(\omega))^4+\frac{16}{3}\omega^3(d(\omega))^3+
\frac{48}{7}\omega^5(d(\omega))^2+4\omega^7d(\omega)+\frac{8}{9}\omega^9).\\
\end{eqnarray*}

\section{(K,N)-flat connections}

Suppose now that we choose a nilpotent derivation
$\delta\in\Omega^1(M,End(\mathcal{E}))$ instead of the exterior
differential, in this case any covariant derivative on $M\times E$
still may be written $\nabla=\delta+\omega$ for some
$\omega\in\Omega^1(M,End(\mathcal{E}))$. Before continue we will
review some notations from [AD].

\smallskip

For $s=(s_1,...,s_n)\in \N^n$ we set $l(s)=n$, the length of the
vector $s$, and $|s|=\sum_i{s_i}$. For $1\leq i <n,\ s_{>i}$
denotes the vector given by $s_{>i}=(s_{i+1},...,s_n)$, for
$1<i\leq n,\ s_{<i}$ stands for $s_{<i}=(s_1,...,s_{i-1})$, we
also set $s_{>n}=s_{<1}=\emptyset$. $\N^{(\infty)}$ denotes the
set $\bigsqcup_{n=0}^{\infty}\N^n$, where by convention
$\N^{0}=\{\emptyset\}$.

\smallskip

For $e\in\Omega^1(M,End(\mathcal{E}))$, we define
$e^{(s)}=e^{(s_1)}...e^{(s_n)}$, where $e^{(a)}=d_{End}^{a}(e)$ if
$a\geq 1$, $e^{(0)}=e$ and $e^{\emptyset}=1$. In the case that
$e_\omega\in\Omega^1(M,End(\mathcal{E}))$ given by
\[e_\omega(\alpha)=\omega\wedge\alpha,\hspace{.3cm}\omega\in\Omega^1(M)\]
then $e_\omega^{(a)}=d_{End}^{a}(e_\omega)$ reduce to
$e_\omega^{(a)}=e_{d^{a}(\omega)}$, thus
\[e_\omega^{(s)}=e_\omega^{(s_1)}\cdots e_\omega^{(s_n)}=e_{d^{s_1}(\omega)}\cdots e_{d^{s_n}(\omega)}.\]

For $N\in\N$ we define $E_N=\{s\in\N^{(\infty)}:|s|+l(s)\leq N\}$
and for $s\in E_N$ we define $N(s)\in\Z$ by $N(s)=N-|s|-l(s)$.

\smallskip

We introduce the discrete quantum mechanical system $L$ by
\begin{enumerate}
\item $V_{L}=\N^{(\infty)}$. \item There is a unique directed edge
in $L$ from vertex $s$ to $t$ if and only if
$t\in\{(0,s),s,(s+e_i)\}$ where
$e_i=(0,..,\underset{\scriptsize{i-th}}{\underbrace{1}},..,0)\in\N^{l(s)}$.

\item Edges in $L$ are weighted according to the following table
\smallskip
\begin{center}
\begin{tabular}{|l|l|l|}\hline $s(e_i)$    &      $t(e_i)$ &
$v(e_i)$      \\  \hline
$s$         &      $(0,s)$      & $1$                 \\
$s$         &      $s$          & $(-1)^{|s|+l(s)}$    \\
$s$         &     $(s+e_i)$    & $(-1)^{|s_{<i}|+i-1}$ \\
\hline
\end{tabular}
\end{center}
\end{enumerate}

The set $P_N(\emptyset,s)$ consist of all paths
$\gamma=(e_1,...,e_N)$, such that $s(e_1)=\emptyset$ and
$t(e_N)=s$. For $\gamma\in P_N(\emptyset,s)$ we define the weight
$v(\gamma)$ of $\gamma$ as
\[v(\gamma)=\prod_{i=1}^{N}v(e_i).\]

For example if we have $\delta^K=0$ then the $N$-curvature is
given by

\begin{theorem}\label{MCMN}
Let $M$ be a manifold and $\mathcal{E}=M\times E$ be a trivial
bundle over $M$. Given any connection $\omega$ consider the
covariant differential $\nabla=\delta+\omega$, then the
N-curvature is given by
\[\nabla^N=\sum_{k=0}^{N-1}c_k\delta^k\]
where
\[c_k=\sum_{\begin{subarray}{c} s\in E_N\\ N(s)=k \\ s_i<K\\ \end{subarray}}c(s,N)\omega^{(s)}
\hspace{.5cm}\text{and}\hspace{.5cm} c(s,N)=\sum_{\gamma\in
P_N(\emptyset,s)}v(\gamma).\]
\end{theorem}

Suppose that we have a derivation $\delta$ such that $\delta^3=0$,
and we want to deform it into a covariant 3-differential
$\nabla=\delta+\omega$, $\omega\in\Omega^1(M)$. So we required
that $\nabla^3=0$. By Theorem \ref{MCMN} we must have
$\displaystyle{\sum_{k=0}^{2}c_k\delta^k=0}$, let us calculate the
coefficients $c_k$. Notice that
\[E_3=\{\emptyset, (0), (1), (2),(0,0), (1,0), (0,1), (0,0,0)\}\]
Let us first compute ${\mathbf{c_0}}$. In this case we have four
vectors in $E_3$ such that $N(s)=0$, these are $(2), (1,0), (0,1),
(0,0,0)$.

\smallskip

For $s=(2)$ there is only one
path from $\emptyset$ to $(2)$ of length 3,\\
$\emptyset\to(0)\to(1)\to(2)$\hspace{.3cm} with weight is 1 and
$\omega^{(2)}=\delta^2(\omega)$, thus $c(s,3)=\delta^2(\omega)$.

\medskip

For $s=(1,0)$ there is only one path from $\emptyset$ to $(1,0)$
of length 3,\\
$\emptyset\to(0)\to(0,0)\to(1,0)$\hspace{.3cm} with weight 1 and
$\omega^{(1,0)}=\delta(\omega)\omega$, thus
$c(s,3)=\delta(\omega)\omega$.

\medskip

For $s=(0,1)$ the paths from $\emptyset$ to $(0,1)$ of length 3
are\\
$\emptyset\to(0)\to(0,0)\to(0,1)$\hspace{.3cm}with weight -1.\\
$\emptyset\to(0)\to(1)\to(0,1)$\hspace{.3cm}with weight 1.\\
$\omega^{(0,1)}=\omega\delta(\omega)$, then $c(s,3)=0$, because
the sum of the weight of the 2 paths is 0.

\medskip

For $s=(0,0,0)$ the only path from $\emptyset$ to $(0,0,0)$ of
length 3 is\\
$\emptyset\to(0)\to(0,0)\to(0,0,0)$\hspace{.3cm}the weight is 1
and $\omega^{(0,0,0)}=\omega^3$, thus $c(s,3)=\omega^3$.

\bigskip

We proceed to compute ${\mathbf{c_1}}$. In this case we have two
vectors in $E_3$ such that $N(s)=1$, $(1)$ and $(0,0)$

\smallskip

For $s=(1)$, the paths from $\emptyset$ to $(1)$ of length 3 are\\
$\emptyset\to\emptyset\to(0)\to(1)$\hspace{.3cm}with weight 1.\\
$\emptyset\to(0)\to(0)\to(1)$\hspace{.3cm}with weight -1.\\
$\emptyset\to(0)\to(1)\to(1)$\hspace{.3cm}with weight 1.\\
$\omega^{(1)}=\delta(\omega)$ and thus $c(s,3)=\delta(\omega)$.

\medskip

For $s=(0,0)$ the paths from $\emptyset$ to $(0,0)$ of length 3
are\\
$\emptyset\to(0)\to(0)\to(0,0)$\hspace{.3cm}with weight -1.\\
$\emptyset\to\emptyset\to(0)\to(0,0)$\hspace{.3cm}with weight
1.\\
$\emptyset\to(0)\to(0,0)\to(0,0)$\hspace{.3cm}with weight 1.\\
$\omega^{(0,0)}=\omega^2$ and thus $c(s,3)=\omega^2$.

\bigskip

Finally we compute ${\mathbf{c_2}}$. In this case we have one
vector in $E_3$ such that $N(s)=2$, $(0)$

\smallskip

For the case $s=(0)$ the paths from $\emptyset$ to $(0)$ of length
3 are\\
$\emptyset\to\emptyset\to\emptyset\to(0)$\hspace{.3cm}with weight 1.\\
$\emptyset\to\emptyset\to(0)\to(0)$\hspace{.3cm}with weight -1.\\
$\emptyset\to(0)\to(0)\to(0)$\hspace{.3cm}with weight 1.\\
$\omega^{(0)}=\omega$ and thus $c(s,3)=\omega$.

\medskip

Now the 3-curvature is given by
\[\nabla^3=(\delta^2(\omega)+\delta(\omega)\omega+\omega^3)+(\delta(\omega)+\omega^2)\delta+\omega \delta^2.\]

\begin{example}
Infinitesimal deformations. Suppose that we have a K-differential
$\delta$ and for a connection $\omega$ we consider the
infinitesimal deformation $\nabla=\delta+t\omega$, where $t^2=0$.
$\omega$ is N-flat if we have $\sum_kc_k\delta^k=0$, but in this
case
\[(t\omega)^{(s)}=(t\omega)^{(s_1)}\cdots(t\omega)^{(s_{l(s)})}=t^{l(s)}\omega^{(s)}=0\hspace{.3cm}\text{unless $l(s)\leq 1$}.\]
Thus $E_N$ is given by
\[E_N=\{(1),\cdots,(K-1)\}\]
and $c_k=c((N-(k+1)),N)\omega^{(N-(k+1))}$ for all $k$ such that
$2\leq N-k\leq K$.
\end{example}

\section{Differential forms of depth N}

We consider $U\subseteq\R^k$ an open set. We use local coordinates
$x_1,\cdots,x_k$ on $U$. For $(N_1,\cdots,N_k)\in\N^k_{\geq 2}$ we
define the algebra $\Omega_{(N_1,\cdots,N_k)}(U)$ of differential
forms of depth $(N_1,\cdots,N_k)$ on $\R^k$ as follows: an element
$\alpha\in \Omega_{(N_1,\cdots,N_k)}(U)$ is given by
\[\alpha=\sum_{I}\alpha_Idx^I\]
where $I:D(I)\subset [k]\to\N^+$, and for all $i\in D(I)$ we have
$I(i)\in[N_i-1]$, $\alpha_I\in C^\infty(U)$ and
$\displaystyle{dx^I=\prod_{i\in D(I)}d^{I(i)}x_i}$.

\smallskip

We declare that $deg(x_i)=0$ and $deg(d^kx_i)=k$ then
$\Omega_{(N_1,\cdots,N_k)}(\R^n)$ is a graded algebra with the
product given by: for
$\alpha,\beta\in\Omega_{(N_1,\cdots,N_k)}(\R^n)$,
$\alpha=\sum_I\alpha_Idx^I$ and $\beta=\sum_J\beta_Jdx^J$
\[\alpha\beta=\sum_{I}(\alpha\beta)_Idx^I,\]
where
\[(\alpha\beta)_I=\sum_{K_1\sqcup K_2=D(I)}sgn(K_1,K_2,I)\alpha_{I|K_1}\beta_{I|K_2},\]
and $\displaystyle{sgn(K_1,K_2,I)=\prod_{\begin{subarray}{c}i\in
K_1,\ j\in K_2
\\ i>j \end{subarray}}(-1)^{I(i)I(j)}}$.

For $\displaystyle{\alpha=\sum_{I}\alpha_{I}dx^I}$ we define the
form $d\alpha$ by
\[d\alpha=\sum_{1\leq s\leq k}\frac{\partial \alpha_I}{\partial x_s}dx_s\wedge dx^I+\sum_{s\in D(I)}(-1)^{\sum_{t<s}I(t)}\alpha_Idx^{I+\delta_s}\]
where $\delta_s:[k]\to\N$ is given by
$e_s(j)=\left\{\begin{array}{cc} 1 & j=s \\ 0 & j\neq s
.\end{array}\right.$

We have the following

\begin{theorem}
$\Omega_{(N_1,\cdots,N_k)}(\R^n)$ is a $(N_1+\cdots+N_k-k)$-dga.
\end{theorem}
{\bf proof} Since
$\Omega_{(N_1,\cdots,N_k)}(x_1,...,x_k)=\Omega_{N_1}(x_1)\hat{\otimes}\cdots\hat{\otimes}\Omega_{N_k}(x_k)$,
we just need to consider the one variable case
$\Omega_N(x)=\R[x,dx,...,d^{N-1}x]/<d^ixd^jx>$, for $1\leq i\leq
N-1$. For $\alpha=\sum_{i=0}^{N-1}f_i(x)d^i(x)$ we have
\[d^k\alpha=\frac{\partial f_0}{\partial x}d^kx+\sum_{i=1}^{N-k-1}f_i(x)d^{i+k}x.\]
Taking $k=N$ we see that $\Omega_N(x)$ is a $N$-complex, the
Leibniz rule and associativity of the product are easy to check,
thus $\Omega_N(x)$ is a N-dga. By [AD, Theorem 8]
$\Omega_{(N_1,\cdots,N_k)}(\R^n)$ is a
$(N_1+\cdots+N_k-k)$-dga.$\blacklozenge$

\smallskip

We shall use the notation
$\Omega_N(\R^k):=\Omega_{N,\cdots,N}(\R^k)$.

Let us recall the definition of affine varieties.

\begin{definition}
$M$ is an affine variety if there is an open covering
$\Lambda=\{U_i\}$ of $M$ and diffeomorphisms
$\varphi_i:U_i\to\R^n$, with $\varphi_i(U_i)$ open, such that
\[\varphi_j\circ\varphi_i^{-1}:\varphi_i(U_i\cap U_j)\longrightarrow \varphi_j(U_i\cap U_j)\]
satisfy $\varphi_j\circ\varphi_i^{-1}(x)=A_{ij}x+b_{ij}$,
$A_{ij}\in Gl_n(\R)$ and $b_{ij}\in\R^n$.
\end{definition}

$f:\R^k\to\R^k$ is an affine map, i.e., $f(x)=Ax+b$ for $A\in
Gl_k(\R^k)$, $b\in\R^k$. The map given by
\[f^*(x_i)=f^i\hspace{.5cm}\text{and}\hspace{.5cm}f^*(d^jx_i)=\frac{\partial f^i}{\partial x_j}d^kx_j,\]
can be extended to a linear map
$f:\Omega_N(\R^k)\to\Omega_N(\R^k)$ such that 1) $df^*=f^*d$ and
2) $f^*(\alpha\wedge\beta)=f^*(\alpha)\wedge f^*(\beta)$, for all
$\alpha,\ \beta\in\Omega_N(\R^k)$.

\smallskip

Let $M$ be an affine manifold. We define the algebra $\Omega_N(M)$
of differential forms of depth N on $M$ as follows, $\Omega_N(M)$
consist of tuples $\alpha=(\alpha_U)_{U\in\Lambda}$, where
$\alpha_U\in\Omega_N(\varphi_U(U))$, satisfying the following
compatibility condition
\[(\varphi_V\circ\varphi_U^{-1})^*(\alpha_V|_{\varphi_V(U\cap
V)})=\alpha_U|_{\varphi_U(U\cap V)}.\] For $U,V\in\Lambda$ such
that $U\cap V\neq\emptyset$. Our final result is the following
\begin{theorem}
Let $M$ an affine manifold with $dim(M)=m$. $\Omega_N(M)$ is a
$(m-1)N$-dga and $\mathbf{H}_N(M)$ is a $((m-1)N-1)$-dga.
\end{theorem}

% ------------------------------------------------------------------------

% ------------------------------------------------------------------------

\[\begin{array}{c}

\mbox{Mauricio Angel. Universidad Central de Venezuela (UCV).} \ \  \mbox{\texttt{mangel@euler.ciens.ucv.ve}} \\

\mbox{Rafael D\'\i az. Universidad Central de Venezuela (UCV).} \ \  \mbox{\texttt{rdiaz@euler.ciens.ucv.ve}} \\
\end{array}\]

\end{sloppypar}

\begin{thebibliography}{00}
\bibitem[A]{} V. Abramov, \textit{On realizations of exterior calculus with dN =0},
Czechoslovak Journal of Physics, 48 (1998) 1265-1272.

\bibitem[AB]{AB} V. Abramov, N. Bazunova, \textit{Exterior calculus with $d^3=0$ on a free associative
algebra and reduced quantum plane}, Proceedings of XIV Max Born
Symposium, New Symmetries and Integrable Models, Karpacz, Poland,
21-24.09.1999, World Scientific, Singapore-New Jersey-London-Hong
Kong, 2000, pp. 3-7.

\bibitem[AK]{AK} V. Abramov, R. Kerner, \textit{On certain realizations of the
q-deformed exterior differential calculus}, Reports on Math.
Phys., 43 (1999) 179-194

\bibitem [AD]{AD} M. Angel, R. D\'{\i}az, \textit{N-differential graded
algebras}, Preprint math.DG/0504398.

\bibitem [BGV]{BGV} N. Berline, E. Getzler, M. Vergne, \textit{Heat Kernels and Dirac
Operators}, Grundlehren 298, Springer-Verlag,
Berlin-Heidelberg-New York, 1992.

\bibitem [F]{F} T. Frankel, \textit{The geometry of physics: an
introduction}, Cambridge University Press, 1997.

\bibitem [DV1]{DV1} M. Dubois-Violette, \textit{Generalized differential spaces with
$d^N=0$ and the q-differential calculus}, Czech J. Phys. 46 (1996)
1227- 1233.

\bibitem[DV2]{DV2} M. Dubois–Violette, \textit{Generalized homologies for $d^N=0$ and
graded q-differential algebras}, Contemporary Mathematics 219,
American Mathematical Society 1998, M. Henneaux, J. Krasil'shchik,
A. Vinogradov, Eds., p. 69-79.

\bibitem[DV3]{DV3} M. Dubois-Violette, \textit{Lectures on differentials,
generalized differentials and some examples related to theoretical
physics}, Contemporary Mathematics 294, American Mathematical
Society 2002, R. Coquereaux, A. Garcia, R. Trinchero Eds, p.
59-94.

\bibitem[DVK]{DVK} M. Dubois-Violette, R. Kerner, \textit{Universal q-differential
calculus and q-analog of homological algebra}, Acta Math. Univ.
Comenianae Vol. LXV, 2 (1996), 175-188.

\bibitem[K]{K} M.M. Kapranov, \textit{On the q-analog of homological
algebra}, Preprint q-alg/9611005.

\bibitem[KW]{KaWa} C. Kassel, M. Wambst, \textit{Alg\`ebre homologique des
N-complexes et homologie de Hochschild aux racines de l'unit\'e},
Publ. Res. Inst. Math. Sci. Kyoto University 34 (1998), nº 2,
91-114.

\bibitem[KN]{KN} R. Kerner, B. Niemeyer, \textit{Covariant q-differential calculus and its
deformations at $q^N=1$}, Lett. in Math. Phys., 45,161-176,
(1998).

\bibitem[M]{Ma} W. Mayer, \textit{A new homology theory I, II},
Annals of Math. 43 (1942) 370-380 and 594-605.

\bibitem[S]{Sit} A. Sitarz; \textit{On the tensor product construction for q-differential
algebras}, Lett. Math. Phys. 44 1998.

\end{thebibliography}
\end{document}